\documentclass[a4paper]{article}

\usepackage{amssymb}
\usepackage{amsthm}
\usepackage{amsmath}

\usepackage{enumerate}
\usepackage{bm}

\newtheorem{theorem*}{Theorem}[section]

\newtheorem{note*}[theorem*]{Note}
\newtheorem{lemma*}[theorem*]{Lemma}
\newtheorem{definition*}[theorem*]{Definition}
\newtheorem{proposition*}[theorem*]{Proposition}
\newtheorem{corollary*}[theorem*]{Corollary} 
\newtheorem{remark*}[theorem*]{Remark}
\newtheorem{example*}[theorem*]{Example}
\numberwithin{equation}{section}

\begin{document}

\title{
Confidence interval for correlation estimator between latent processes
}

\author{
Akitoshi Kimura 
\thanks{Graduate School of Mathematical Sciences, University of Tokyo, 3-8-1 Komaba, Meguro-ku, Tokyo 153-8914, Japan, 
Japan Science and Technology Agency, CREST}
}
%
%

\maketitle

\begin{abstract}
Kimura and Yoshida \cite{Kimura2016} treated a model 
in which 
the finite variation part of a two-dimensional semimartingale is expressed by time-integration of latent 
processes. 
%
They proposed a correlation estimator between the latent processes and proved its consistency and asymptotic mixed normality.  
In this paper, we discuss the confidence interval of the correlation estimator to detect the correlation. 
We propose two types of estimators for asymptotic variance of the correlation estimator and prove their consistency in a high frequency setting. 
Our model includes 
doubly stochastic Poisson processes whose intensity processes are correlated It\^o processes. 
We compare 
our estimators based on 
the simulation of the 
doubly stochastic Poisson processes.  
%
\end{abstract}

\section{Introduction}
%
%
Integrated correlation is an important index in high frequency financial data analysis.
Epps \cite{epps1979comovements} pointed out that the sample correlation between the returns of two different stocks decreases as the sampling frequency of data increases. 
It is considered that non-synchronicity and market microstructure of trading cause this phenomenon. 

Non-synchronous covariance estimation schemes 
have been developed: 
Fourier analytic approach 
(Malliavin and Mancino \cite{malliavin2002fourier}, 
Malliavin et al. \cite{Malliavin2007}) and 
the cumulative covariance estimator 
(Hayashi and Yoshida \cite{
hayashi2005covariance,hayashi2008asymptotic,hayashi2011nonsynchronous}, 
Mykland \cite{mykland2012gaussian}). 

The 
market microstructure is 
modeled 
as the noise added to the latent price process. 
This modeling was successful and
denoising techniques have also been developed: 
sub-sampling 
(Zhang et al. \cite{
zhang2012tale}, 
Zhang \cite{zhang2006efficient}), 
pre-averaging 
(Podolskij and Vetter \cite{
podolskij2007estimation},
Jacod et al. \cite{jacod2009microstructure}), 
and others (Zhou \cite{zhou1996high}). 
There are many studies that treat both non-synchoronicity and market microstructure noise: 
Malliavin and Mancino \cite{malliavin2009fourier}, 
Mancino and Sanfelici \cite{
mancino2011estimating}, 
Park and Linton \cite{park2012estimating}, 
Voev and Lunde \cite{voev2007integrated}, 
Griffin and Oomen \cite{griffin2011covariance}, 
Christensen et al. \cite{christensen2010pre, 
christensen2013covariation},  
Koike \cite{koike2013estimation, koike2014limit,koike2014estimator},  
A{\"\i}t-Sahalia et al. \cite{ait2010high}, 
Barndorff-Nielsen et al. \cite{barndorff2011multivariate},  
Bibinger \cite{bibinger2011efficient,bibinger2012estimator}. 
The followings are proposed as the cause of market microstructure: 
bid-ask spread
(Roll \cite{10.2307/2327617}), 
discretization error
(Gottlieb and Kalay \cite{10.2307/2328052}),
and 
asymmetric information
(Glosten and Milgrom \cite{GLOSTEN198571}). 
On the other hand, 
the relationship between the additive noise modeling and the 
market microstructure is not so clearly explained. 

In ultra high frequency sampling, recently, 
the market microstructure is modeled as dynamics of the limit order book (LOB) rather than the noise. 
Many studies adopt approaches modeling LOB with 
Poisson processes 
(Cont et al. \cite{cont2010stochastic}, 
Abergel and Jedidi \cite{abergel2013mathematical}, 
Muni Toke 
\cite{doi:10.1080/14697688.2014.963654, 2015arXiv150203871M}, 
Smith et al. \cite{smith2003statistical}, 
Muni Toke and Yoshida 
\cite{2016arXiv160203944M}), 
Hawkes processes
(Hewlett \cite{hewlett2006clustering}, 
Large \cite{large2007measuring}, 
Bowsher \cite{bowsher2007modelling}, 
Bacry et al. \cite{bacry2013modelling}, 
Abergel and Jedidi \cite{
abergel2015long}, 
Muni Toke and Pomponio \cite{muni2011modelling}, 
Muni Toke 
\cite{
2010arXiv1003.3796M}, 
Clinet and Yoshida 
\cite{CLINET20171800}, 
Ogihara and Yoshida \cite{
ogihara2015quasi}), and
doubly stochastic Poisson processes 
(Abergel et al. \cite{2017arXiv170501446A}, 
Guilbaud and Pham \cite{2011arXiv1106.5040G}, 
Chertok et al. \cite{Chertok2016}, 
Korolev et al. \cite{2014arXiv1410.1900K}).  
These approaches can also treat the non-synchronicity of trading. 

In this stream, Kimura and Yoshida \cite{Kimura2016} focused on the integrated correlation between intensity processes of doubly stochastic Poisson processes. 
They introduced the following model to treat it in a generalized form. 
This paper is also based on the model. 

Now, we consider a stochastic basis 
$\mathcal{B} = (\Omega, \mathcal{F}, \mathbb{F}, P)$, 
$\mathbb{F} = (\mathcal{F}_{t})_{t \in [0,T]}$. 
On $\mathcal{B}$, let
$\mathbb{X} = (X^1, X^2)$ be an $\mathbb{R}^2$-valued It\^o process given by 
\begin{align} \label{Xint}
\mathbb{X}_{t} = \mathbb{X}_{0} + \int_{0}^{t} \mathbb{X}^{0}_{s} \ ds + \int_{0}^{t} \mathbb{X}^{1}_{s} \ dw_{s}
\ \ (t \in [0,T]), 
\end{align}
where 
$w$ is an $r$-dimensional $\mathcal{F}$-Wiener process, 
$\mathbb{X}_{0}$ is an $\mathbb{F}_{0}$-measurable random variable, 
$\mathbb{X}^{0}$ is a two-dimensional $\mathbb{F}$-adapted process, and 
$\mathbb{X}^{1}$ is an $\mathbb{R}^{2} \otimes \mathbb{R}^{r}$-valued $\mathbb{F}$-adapted process 
satisfies condition [A] mentioned later. 
Let
$a_n$ be a positive number depending on $n$. 
On $\mathcal{B}$, consider a two-dimensional measurable process 
$\mathbb{Y}^n = (Y^{n, 1}, Y^{n, 2})$ having a decomposition 
\begin{align*}
\mathbb{Y}^{n}_{t} = \mathbb{Y}^{n}_{0} + \int_{0}^{t} a_n \mathbb{X}_{s} \ ds + \mathbb{M}^{n}_{t} 
\ \ (t \in [0,T]),   
\end{align*}
where $\mathbb{M}^{n}$ is a two-dimensional measurable process with $\mathbb{M}^{n}_{0} = 0$ 
satisfies condition 
[B], [B'], or [B$^{\sharp^{\sharp}}$] mentioned later.  

In this model,  
$\mathbb{Y}$ can be a counting process with intensity process $a_n \mathbb{X}$ 
where $\mathbb{X}$ is $\mathbb{R}_+^2$-valued, 
i.e. $\mathbb{Y}$ is a two-dimensional doubly stochastic Poisson process. 
This counting process models the high frequency counting data of the orders or transactions in the active market, for example. 

Let $I_j = [t_{j-1}, t_{j})$ for a sampling design $\Pi = (t_j)_{j=0, \dots , b_n}$ 
with $0 = t_{0} < t_{1} < \dots  < t_{b_n} =T$ 
and $h_{j} = t_{j} - t_{j-1}$.  

Kimura and Yoshida \cite{Kimura2016} introduced the 
estimator of the covariance $\langle X^1, X^2 \rangle$ from the sampled data of $\mathbb{Y}$ 
as follows.  
For $\alpha, \beta = 1, 2$, 
\begin{align*}
S_n^{\alpha \beta} 
&= \sum_{k=2}^{b_n} 
\Big(\frac{Y^{n, \alpha}_{t_{k}} - Y^{n, \alpha}_{t_{k-1}}}{a_n h_{k}} - \frac{Y^{n, \alpha}_{t_{k-1}} - Y^{n, \alpha}_{t_{k-2}}}{a_n h_{k-1}} \Big) 
\Big(\frac{Y^{n, \beta}_{t_{k}} - Y^{n, \beta}_{t_{k-1}}}{a_n h_{k}} - \frac{Y^{n, \beta}_{t_{k-1}} - Y^{n, \beta}_{t_{k-2}}}{a_n h_{k-1}} \Big), 
\end{align*}
Here, $S_n = (S_n^{1 2}, S_n^{1 1}, S_n^{2 2})^{\star}$ is the (co)variance estimator. 
$S_n$ depends on the scaling parameter $a_n$ which is not derived from the data of $\mathbb{Y}$. 
Therefore, treat the correlation estimator $C_n^{1 2} = S_n^{1 2}/\sqrt{S_n^{1 1} S_n^{2 2}}$ 
which does not depend on the scaling parameter $a_n$. 

Let
\begin{align*}
U^{\alpha \beta} = \frac{2}{3} \langle X^{\alpha,c}, X^{\beta,c} \rangle_T
\>=\>\frac{2}{3}\int_0^T X^{\alpha1}_{t}\cdot X^{\beta1}_{t}dt
\quad(\alpha,\beta=1,2)
\end{align*}
and  
$U = (U^{12}, U^{11}, U^{22})^\star$, 
where $X^{\alpha,c}$ is the continuous part of $X^{\alpha}$ 
and $X^{\alpha1}_{t}$ is the $\alpha$-th row of $\mathbb{X}_t^{1}$.
%
%
%
Let 
\begin{align*} 
\gamma^{(\alpha_1,\beta_1),(\alpha_2,\beta_2)} 
&= \int_0^T \sum_{i,j=1}^{r} 
\frac{X_{i,s}^{\alpha_1 1} X_{j,s}^{\beta_1 1} + X_{j,s}^{\alpha_1 1} X_{i,s}^{\beta_1 1}}{2} 
\frac{X_{i,s}^{\alpha_2 1} X_{j,s}^{\beta_2 1} + X_{j,s}^{\alpha_2 1} X_{i,s}^{\beta_2 1}}{2} \ ds \\
&= \int_0^T (\mathbb{X}_s^{1})^{\tilde{\otimes} (\alpha_1,\beta_1)} 
\cdot (\mathbb{X}_s^{1})^{\tilde{\otimes} (\alpha_2,\beta_2)} \ ds.
\end{align*}
and $\Gamma = (\gamma ^{p q})_{p,q = (1, 2), (1,1), (2,2)}$, 
where $X_{i,s}^{\alpha 1}$ is the $(\alpha, i)$-element of  $\mathbb{X}_s^{1}$. 


Under [A] and any one of [B], [B'], and [B$^{\sharp^{\sharp}}$], 
Kimura and Yoshida \cite{Kimura2016} proved the the asymptotic mixed normality of the covariance estimator and the correlation estimator. 
\begin{align*}
\Big(\frac{T}{b_n} \Big)^{-1/2} (S_n - U) \to^{d_s} \Gamma^{1/2} \zeta
\end{align*}
as $n \to \infty$,  
where $\zeta$ is an $\mathbb{R}^3$-valued standard normal variable independent of $\mathcal{F}$. 
and $d_s$ means $\mathcal{F}$-stable convergence. 
Let $R = U^{12} / \sqrt{U^{11} U^{22}}$
and suppose that $U^{11} U^{22} \neq 0$ a.s.  
\begin{align*}
\Big(\frac{T}{b_n} \Big)^{-1/2} (C_n^{12} - R) \to^{d_s} \Xi^{1/2} \zeta
\end{align*}
as $n \to \infty$, where 
\begin{align*}
\Xi^{1/2} := 
\Big(\frac{1}{\sqrt{U^{11} U^{22}}}, \frac{-U^{12}}{2 \sqrt{(U^{11})^3 U^{22}}}, \frac{-U^{12}}{2 \sqrt{U^{11} (U^{22})^3}} \Big) \Gamma^{1/2}. 
\end{align*}
$\Xi = (\Xi^{1/2}) (\Xi^{1/2})^{\star}$ is the asymptotic variance of the correlation estimator. 

In this paper, we introduce two types of estimators for $\Gamma$ and prove their consistency 
in order to obtain the consistent estimator for  $\Xi$. 
Therefore, we can detect the correlation between latent processes based on hypothesis testing. 

This paper is organized as follows.
In Section \ref{assumption}, we list the assumptions. 
In Section \ref{result}, we introduce the 
estimators for $\Gamma$ and $\Xi$, and state the main result: their consistency. 
In Section \ref{proof}, the results in Section \ref{result} are proved.  
In Section \ref{simulation}, an example and simulation studies are given.  

\section{Assumptions} \label{assumption}
For simplicity,  we assume 
$t_{j} = j (T/b_n)$, $h_{j} = T/b_n =: \delta_n$. 
We write $\Delta_{j} V = V_{t_{j}} - V_{t_{j-1}}$ for a process $V$. 
Consider the following conditions.  
\begin{itemize}
\item[\bf{[A]}]
Process $\mathbb{X}$ admits the representation (\ref{Xint}) 
for an $\mathbb{R}^2$-valued $\mathcal{F}_0$-measurable random variable $\mathbb{X}_0$ 
and coefficients $\mathbb{X}^{\kappa}$ ($\kappa=0,1$) such that 
$\mathbb{X}^{0}$ is a c\'adl\'ag $\mathbb{F}$-adapted process and that 
$\mathbb{X}^{1}$ has a representation
\begin{align*}
\mathbb{X}^{1}_{t} 
= \mathbb{X}^{1}_{0} + \int_{0}^{t} \mathbb{X}^{10}_{s} \ ds + \int_{0}^{t} \sum_{\kappa' = 1}^{r'} \mathbb{X}^{1\kappa'}_{s} \ d\tilde{w}^{\kappa'}_s 
\ \ (t \in [0,T]), 
\end{align*}
where
$\mathbb{X}^{1}_{0}$ is an $\mathbb{R}^{2} \otimes \mathbb{R}^{r}$-valued $\mathcal{F}_{0}$-measurable random variable, 
$\tilde{w} = (\tilde{w}^1,...,\tilde{w}^{r'})$ is an $r'$-dimensional $\mathbb{F}$-Wiener process (not necessary independent of $w$), 
and
$\mathbb{X}^{1 \kappa'}_{s}$ $(\kappa' = 0, 1, ..., r')$ are $\mathbb{R}^{2} \otimes \mathbb{R}^{r}$-valued c\'adl\'ag $\mathbb{F}$-adapted processes. 

%
%
\item[\bf{
[B]
}]
$\mathbb{M}^n = (M^{n, \alpha})_{\alpha=1, 2}$ is a two-dimensional $\mathbb{F}$-local martingale 
with $\mathbb{M}^n_{0} = 0$ and such that  
\begin{itemize}
\item[
(i)
] 
$\lim_{n \to \infty} b_n^{5/2} / a_n = 0$. 

\item[
(ii)
] 
$\sum_{j=1}^{b_n} |\Delta_{j} \mathbb{M}^n|^2 = O_p(a_n)$ as $n \to \infty$, 
$\sup_{t \in [0,1]}  |\Delta \mathbb{M}^n| \leq c a_n^{1/2}$ for a constant c independent of n. 

\item[
(iii)
] 
The absolutely continuous (w.r.t. the Lebesgue measure a.s.) mapping 
$[0, T] \ni t \mapsto \langle \mathbb{M}^n , w \rangle_{t} \in \mathbb{R}^{2} \otimes \mathbb{R}^{r}$ 
satisfies 
$\sup_{t \in [0, T]} |d \langle \mathbb{M}^n , w \rangle_{t} / dt| = O_p(b_n)$ as $n \to \infty$. 
\end{itemize}

Here, 
$\langle \mathbb{M}^n , w \rangle_{t}$ is the $2 \times r$ matrix of angle brackets $\langle M^{n, \alpha} , w^{k} \rangle_{t}$ for $\mathbb{M}^n = (M^{n, \alpha})_{\alpha = 1, 2}$ and $w = (w^k)_{k=1,...,r}$. 

\item[\bf{
[B$'$]
}]
$\mathbb{M}^n$ is a two-dimensional $\mathbb{F}$-local martingale 
with $\mathbb{M}^n_{0} = 0$, satisfies [B] (i), (iii), and 
\begin{itemize}
\item[
(ii$'$)
]
$E[\sum_{j=1}^{b_n} |\Delta_{j} \mathbb{M}^n|^2] = O(a_n)$ as $n \to \infty$. 
\end{itemize}

\item[\bf{
[B$^{\sharp^{\sharp}}$]
}]
$\mathbb{M}^n$ is a two-dimensional $\mathbb{F}$-local martingale 
with $\mathbb{M}^n_{0} = 0$, satisfies [B] (ii), (iii), and 
\begin{itemize}
\item[
(i$^{\sharp^{\sharp}}$)
]
$\lim_{n \to \infty} b_n^{3} / a_n = 0$. 
%
\end{itemize}

\item[\bf{
[C]
}]
$E[\sum_{j=1}^{b_n} |\Delta_{j} \mathbb{M}^n|^4] = O(a_n^2 b_n^{-1})$ as $n \to \infty$
\end{itemize}

\begin{remark*}

[A], [B], and [B$'$] are the same condition in \cite{Kimura2016}. 
[B$^{\sharp^{\sharp}}$] is a bit stronger condition than [B$^{\sharp}$] in \cite{Kimura2016}. 
\end{remark*}

\begin{remark*}

The doubly stochastic Poisson process model satisfies [C]. 
\end{remark*}

\section{Results} \label{result}
We need the estimator of $\Gamma$ to 
obtain the asymptotic variance estimator. 
Here, we give two types of the gamma estimator. 

Before that, we introduce some notations to simplify the description. 
Let 
$A^{\otimes(i, j)}$ is an $(i, j)$-element of $A^{\otimes}$, 
$A^{\otimes} = A \otimes A = AA^{\star}$, and 
$^\star$ denotes transpose. 
Let $x\tilde{\otimes}y=((x_iy_j+x_jy_i)/2)\in\mathbb{R}^{\sf r}\otimes\mathbb{R}^{\sf r}$ 
for $x=(x_i)$, $y=(y_i)\in\mathbb{R}^{\sf r}$,  and 
let $x^{\tilde{\otimes}(\alpha,\beta)}=x^\alpha_\cdot\tilde{\otimes}x^\beta_\cdot$ for $x=(x^\alpha_i)\in\mathbb{R}^2\otimes\mathbb{R}^{\sf r}$. 
We write $x\cdot y=\sum_{i=1}^{\sf r}x_iy_i$ for $x=(x_i)$, $y=(y_i)\in\mathbb{R}^{\sf r}$, 
and $x\cdot y = \sum_{i,j=1}^{\sf r}x_{i,j}y_{i,j}$ for 
$x=(x_{i,j})$, $y=(y_{i,j})\in\mathbb{R}^{\sf r}\otimes \mathbb{R}^{\sf r}$. 
Write 
\begin{align*}
\tilde{V}_{k} := \frac{\Delta_{k}V}{a_n h_{k}} = \frac{\Delta_{k}V}{a_n \delta_n}. 
\end{align*}
for a stochastic process $V$. 
For example, we can write 
\begin{align*}
S_n^{\alpha \beta} 
= \sum_{j=2}^{b_n} (\tilde{\mathbb{Y}}^n_{j} - \tilde{\mathbb{Y}}^n_{j-1})^{\otimes(\alpha, \beta)}, 
\quad
\gamma ^{p q} = \int_0^1 (\mathbb{X}_s^{1})^{\tilde{\otimes} p} \cdot 
(\mathbb{X}_s^{1})^{\tilde{\otimes}q} \ ds. 
\end{align*}

\subsection{$\Gamma$ estimator}
In the similar way in 
Barndorff-Nielsen and Shephard \cite{10.2307/3598838}, 
we define 
\begin{align*}
\hat{\gamma}_{n, 1}^{pq} 
&=\frac{9}{8} \Big\{ \sum_{k=2}^{b_n} 
(\tilde{\mathbb{Y}}^n_{k} - \tilde{\mathbb{Y}}^n_{k-1})^{\otimes p} 
(\tilde{\mathbb{Y}}^n_{k} - \tilde{\mathbb{Y}}^n_{k-1})^{\otimes q} \\
&\qquad
- \frac{1}{2} \sum_{k=2}^{b_n - 2}
\Big( (\tilde{\mathbb{Y}}^n_{k} - \tilde{\mathbb{Y}}^n_{k-1})^{\otimes p} 
(\tilde{\mathbb{Y}}^n_{k+2} - \tilde{\mathbb{Y}}^n_{k+1})^{\otimes q} \\ 
&\qquad\quad 
+ (\tilde{\mathbb{Y}}^n_{k+2} - \tilde{\mathbb{Y}}^n_{k+1})^{\otimes p} 
(\tilde{\mathbb{Y}}^n_{k} - \tilde{\mathbb{Y}}^n_{k-1})^{\otimes q} \Big) \Big\} \Big(\frac{T}{b_n} \Big)^{-1}
\end{align*}
By the similar idea, we define 
\begin{align*}
\hat{\gamma}_{n, 2}^{pq} 
&= \frac{9}{8} \sum_{k=2}^{b_n-2} \frac{1}{2} 
\Big\{ (\tilde{\mathbb{Y}}^n_{k+2} - \tilde{\mathbb{Y}}^n_{k+1})^{\otimes p} 
- (\tilde{\mathbb{Y}}^n_{k} - \tilde{\mathbb{Y}}^n_{k-1})^{\otimes p} \Big\} \\
&\qquad\qquad
\times \Big\{ (\tilde{\mathbb{Y}}^n_{k+2} - \tilde{\mathbb{Y}}^n_{k+1})^{\otimes q} 
- (\tilde{\mathbb{Y}}^n_{k} - \tilde{\mathbb{Y}}^n_{k-1})^{\otimes q} \Big\} \Big(\frac{T}{b_n} \Big)^{-1}.  
\end{align*}
%
%

\subsection{Kernel based $\Gamma$ estimator}
In the similar way in 
Hayashi and Yoshida \cite{hayashi2011nonsynchronous}, we define the kernel based estimator. 
\begin{align*}
\partial_{h} \{X^{\alpha} , X^{\beta} \}_{k} 
:= \sum_{l=(k - n(h) +1) \vee 2}^{k} (\Delta_{l} X^{\alpha} \Delta_{l} X^{\beta}) h^{-1}. 
\end{align*}
where $h = h_n$ is a parameter satisfying $h_n \to 0$ and $h_n^{-1} b_n^{-1} \to 0$, 
and $n(h) := \max_{m} \{ t_m \leq h\} $. 
%
%
%
\begin{align*}
\hat{\gamma}_{n, h}^{pq} 
=& \frac{9}{8}\sum_{k=2}^{b_n} \Big\{ 
\partial_{h} \{ \tilde{Y}^{n, \alpha_1} , \tilde{Y}^{n, \alpha_2} \}_{k} \
\partial_{h} \{ \tilde{Y}^{n, \beta_1} , \tilde{Y}^{n, \beta_2} \}_{k} \\
&\qquad
+ \partial_{h} \{ \tilde{Y}^{n, \alpha_1} , \tilde{Y}^{n, \beta_2} \}_{k} \
\partial_{h} \{ \tilde{Y}^{n, \beta_1} , \tilde{Y}^{n, \alpha_2} \}_{k} \Big\} \Big(\frac{T}{b_n} \Big).
\end{align*} 

\subsection{$\Xi$ estimator}
Here, we define 
\begin{align*}
\hat{\Xi}_{n, *}^{1/2} 
= \Big(\frac{1}{\sqrt{S_n^{11} S_n^{22}}}, \frac{-S_n^{12}}{2 \sqrt{(S_n^{11})^3 S_n^{22}}}, \frac{-S_n^{12}}{2 \sqrt{S_n^{11} (S_n^{22})^3}} \Big) \hat{\Gamma}_{n, *}^{1/2} 
\end{align*}
and 
$\hat{\Xi}_{n, *} = (\hat{\Xi}_{n, *}^{1/2})^{\otimes}$, for $* = 1, 2, h$. 

\subsection{consistency of the estimators}
\begin{theorem*} \label{nonkernel}
(a)
$
\hat{\gamma}_{n, i}^{p q} = \gamma^{p q} + o_p(1)$ 
under [A] and any one of [B], [B$'$] and [B$^{\sharp^{\sharp}}$], ($i=1, 2$). 

(b)
$
\hat{\gamma}_{n, i}^{p q} = \gamma^{p q} + O_p(b_n^{-1/2})$ 
under [A] and [B$^{\sharp^{\sharp}}$], ($i=1, 2$). 

(c)
$
\hat{\gamma}_{n, i}^{p q} = \gamma^{p q} + O_p(b_n^{-1/2})$ 
under [A], [B] or [B'], and [C] (i=1,2).  
\end{theorem*}
%
%
\begin{theorem*} \label{kernel}
(a)
$
\hat{\gamma}_{n, h}^{p q} = \gamma^{p q} + o_p(1)$ 
under [A] and any one of [B], [B$'$] and [B$^{\sharp^{\sharp}}$]. 

(b)
$
\hat{\gamma}_{n, h}^{p q} = \gamma^{p q} + O_p(b_n^{-1/2}) + O_p(h_n)$ 
under 
[A] and [B$^{\sharp^{\sharp}}$]. 

(c)
$
\hat{\gamma}_{n, h}^{p q} = \gamma^{p q} + O_p(b_n^{-1/2}) + O_p(h_n)$ 
under 
[A], [B] or [B'], and [C]. 
\end{theorem*}

\begin{corollary*} \label{avar}
(a)
$
\hat{\Xi}_{n, i} = \Xi + o_p(1)$ 
under [A] and any one of [B], [B$'$] and [B$^{\sharp^{\sharp}}$], ($i=1, 2$). 

(b)
$
\hat{\Xi}_{n, i} = \Xi + O_p(b_n^{-1/2})$ 
under [A] and [B$^{\sharp^{\sharp}}$], ($i=1, 2$). 

(c)
$
\hat{\Xi}_{n, i} = \Xi + O_p(b_n^{-1/2})$ 
under [A], [B] or [B'], and [C] (i=1,2).  
\end{corollary*}

\begin{corollary*} \label{avar_k}
(a)
$
\hat{\Xi}_{n, h} = \Xi + o_p(1)$ 
under [A] and any one of [B], [B$'$] and [B$^{\sharp^{\sharp}}$]. 

(b)
$
\hat{\Xi}_{n, h} = \Xi + O_p(b_n^{-1/2}) + O_p(h_n)$ 
under 
[A] and [B$^{\sharp^{\sharp}}$]. 

(c)
$
\hat{\Xi}_{n, h} = \Xi + O_p(b_n^{-1/2}) + O_p(h_n)$ 
under 
[A], [B] or [B'], and [C]. 
\end{corollary*}

\section{Proof} \label{proof}

\subsection*{Proof of Theorem \ref{nonkernel}}
First, we approximate $\hat{\gamma}_{n, i}^{p q} (i = 1, 2)$ by the following quantities. 
Let
\begin{align*}
G_{n,1}^{pq} 
&= \frac{9}{8} \Big\{ \sum_{k=2}^{b_n} 
(\tilde{\bm{\chi}}^n_{k} - \tilde{\bm{\chi}}^n_{k-1})^{\otimes p} 
(\tilde{\bm{\chi}}^n_{k} - \tilde{\bm{\chi}}^n_{k-1})^{\otimes q} \\
&\qquad
- \frac{1}{2} \sum_{k=2}^{b_n - 2}
\Big( (\tilde{\bm{\chi}}^n_{k} - \tilde{\bm{\chi}}^n_{k-1})^{\otimes p} 
(\tilde{\bm{\chi}}^n_{k+2} - \tilde{\bm{\chi}}^n_{k+1})^{\otimes q} \\
&\qquad\qquad\qquad 
+ (\tilde{\bm{\chi}}^n_{k+2} - \tilde{\bm{\chi}}^n_{k+1})^{\otimes p} 
(\tilde{\bm{\chi}}^n_{k} - \tilde{\bm{\chi}}^n_{k-1})^{\otimes q} \Big) \Big\} \Big(\frac{T}{b_n} \Big)^{-1}
\end{align*}
and
\begin{align*}
G_{n,2}^{pq} 
&= \frac{9}{8} \sum_{k=2}^{b_n-2} \frac{1}{2} 
\Big\{ (\tilde{\bm{\chi}}^n_{k+2} - \tilde{\bm{\chi}}^n_{k+1})^{\otimes p} 
- (\tilde{\bm{\chi}}^n_{k} - \tilde{\bm{\chi}}^n_{k-1})^{\otimes p} \Big\} \\
&\qquad\qquad
\times
\Big\{ (\tilde{\bm{\chi}}^n_{k+2} - \tilde{\bm{\chi}}^n_{k+1})^{\otimes q} 
- (\tilde{\bm{\chi}}^n_{k} - \tilde{\bm{\chi}}^n_{k-1})^{\otimes q} \Big\} \Big(\frac{T}{b_n} \Big)^{-1}. 
\end{align*}
where 
$\bm{\chi}^n = (\chi^{n, 1}, \chi^{n, 2})$ and 
$\chi^{n, \alpha}_t = \int_{0}^{t} a_n X^{\alpha}_{s} \ ds$. 
%
%

\begin{lemma*} \label{lemma1i}
(a)
$
\hat{\gamma}^{p q}_{n,i} = 
G_{n,i}^{p q} + o_p(1)$ 
under any one of [B], [B$'$] and [B$^{\sharp^{\sharp}}$], (i=1,2).  

(b)
$
\hat{\gamma}^{p q}_{n,i} = 
G_{n,i}^{p q} + o_p(b_n^{-1/2})$ 
under [B$^{\sharp^{\sharp}}$], (i=1,2).  
\end{lemma*}

\proof
Write
$(V^1 \otimes V^2 \otimes \cdots \otimes V^n)_{i_1, i_2, \ldots, i_n} = V^1_{i_1} V^2_{i_2} \cdots V^n_{i_n}$
for vectors $V^1, V^2, \ldots, V^n$ 
and confirm 
$|V^1 \otimes V^2 \otimes \cdots \otimes V^n | = |V^1| |V^2| \cdots |V^n|$. 
Let
$s_{k} = (\tilde{\mathbb{Y}}^n_{k} - \tilde{\mathbb{Y}}^n_{k-1})^{\otimes}$, 
$t_{k} = (\tilde{\bm{\chi}}^n_{k} - \tilde{\bm{\chi}}^n_{k-1})^{\otimes}$,
%
%
%
$\hat{\gamma}_{n,i} 
= (\hat{\gamma}^{\alpha_1, \beta_1, \alpha_2, \beta_2}_{n,i})_{\alpha_1, \beta_1, \alpha_2, \beta_2 = 1, 2}$, 
and 
${G}_{n,i} 
= ({G}^{\alpha_1, \beta_1, \alpha_2, \beta_2}_{n,i})_{\alpha_1, \beta_1, \alpha_2, \beta_2 = 1, 2}$ ($i=1, 2$). 
%
It holds that
\begin{align*}
&\Big(\frac{9}{8} \Big)^{-1} 
| \hat{\gamma}_{n,1} - G_{n,1} | \\
&= \Big(\frac{T}{b_n} \Big)^{-1} \Big| \sum_{k=2}^{b_n} (s_{k} \otimes s_{k} - t_{k} \otimes t_{k}) \\
&\qquad\qquad\quad
- \frac{1}{2} \sum_{k=2}^{b_n - 2} \Big( (s_{k} \otimes s_{k+2} - t_{k} \otimes t_{k+2}) 
+  (s_{k+2} \otimes s_{k} - t_{k+2} \otimes t_{k}) \Big) \Big| \\
&\leq \Big(\frac{T}{b_n} \Big)^{-1} \Big\{ \sum_{k=2}^{b_n} |s_{k} \otimes s_{k} - t_{k} \otimes t_{k}| \\
&\qquad\qquad\quad
+ \frac{1}{2} \sum_{k=2}^{b_n - 2} |s_{k} \otimes s_{k+2} - t_{k} \otimes t_{k+2}| 
+ \frac{1}{2} \sum_{k=2}^{b_n - 2} |s_{k+2} \otimes s_{k} - t_{k+2} \otimes t_{k}| \Big\} 
\end{align*}
and
\begin{align*}
&\Big(\frac{9}{8} \Big)^{-1} 
| \hat{\gamma}_{n,2} - G_{n,2} | \\
&= \frac{1}{2} \Big(\frac{T}{b_n} \Big)^{-1} \Big| \sum_{k=2}^{b_n-2} 
(s_{k+2} \otimes s_{k+2} - t_{k+2} \otimes t_{k+2}) 
- (s_{k+2} \otimes s_{k} - t_{k+2} \otimes t_{k}) \\
&\qquad\qquad\qquad\quad
- (s_{k} \otimes s_{k+2} - t_{k} \otimes t_{k+2}) 
+ (s_{k} \otimes s_{k} - t_{k} \otimes t_{k}) \Big| \\
&\leq \frac{1}{2} \Big(\frac{T}{b_n} \Big)^{-1} \Big\{ \sum_{k=2}^{b_n-2} 
|s_{k+2} \otimes s_{k+2} - t_{k+2} \otimes t_{k+2}| 
+ \sum_{k=2}^{b_n-2} |s_{k+2} \otimes s_{k} - t_{k+2} \otimes t_{k}| \\
&\qquad\qquad\qquad
+ \sum_{k=2}^{b_n-2} |s_{k} \otimes s_{k+2} - t_{k} \otimes t_{k+2}| 
+ \sum_{k=2}^{b_n-2} |s_{k} \otimes s_{k} - t_{k} \otimes t_{k}|  \Big\}. 
\end{align*}
Now, let $\rho=0, 1/2$ and see the typical term. 
\begin{align*}
&b_n^{1+\rho} \sum_{k=2}^{b_n-2} 
|s_{k+2} \otimes s_{k} - t_{k+2} \otimes t_{k}| \\
&= b_n^{1+\rho} \sum_{k=2}^{b_n-2} 
\Big|(s_{k+2} - t_{k+2}) \otimes (s_{k} - t_{k}) \\
&\qquad\qquad\quad
+ (s_{k+2} - t_{k+2}) \otimes t_{k}
+ t_{k+2} \otimes (s_{k} - t_{k}) \Big| \\
&\leq b_n^{1+\rho} \sum_{k=2}^{b_n-2} 
\Big|(s_{k+2} - t_{k+2}) \otimes (s_{k} - t_{k}) \Big| \\
&\qquad\qquad\quad
+ \Big|(s_{k+2} - t_{k+2}) \otimes t_{k} \Big|
+ \Big|t_{k+2} \otimes (s_{k} - t_{k}) \Big| \\
&=  b_n^{1+\rho} \sum_{k=2}^{b_n-2} 
|s_{k+2} - t_{k+2}| |s_{k} - t_{k}|
+ |s_{k+2} - t_{k+2}| |t_{k}|
+ |t_{k+2}| |s_{k} - t_{k}| \\
&\leq \Big(b_n^{1/2 + \rho/2} \sum_{k=2}^{b_n-2} |s_{k+2} - t_{k+2}| \Big) 
\Big(b_n^{1/2 + \rho/2} \sum_{k=2}^{b_n-2} |s_{k} - t_{k}| \Big) \\
&\quad
+ \Big(b_n^{\rho} \sum_{k=2}^{b_n-2} |s_{k+2} - t_{k+2}|\Big) \Big(b_n \max_{k} |t_{k}| \Big) \\
&\quad
+ \Big(b_n \max_{k} |t_{k+2}| \Big) \Big(b_n^{\rho} \sum_{k=2}^{b_n-2} |s_{k} - t_{k}| \Big).
%
%
\end{align*}
Under [B$^{\sharp^{\sharp}}$], 
the same argument of the proof of Lemma 1 (b) in \cite{Kimura2016} under the assumption [B] yields  
\begin{align*}
b_n^{3/4} \sum_{k=2}^{b_n} |s_{k} - t_{k}| \to^p 0.  
\end{align*}
Under any one of [B], [B$'$] and [B$^{\sharp^{\sharp}}$], it holds that 
\begin{align*}
b_n^{1/2} \sum_{k=2}^{b_n} |s_{k} - t_{k}| \to^p 0 
\end{align*}
by the proof of Lemma 1 (b) in \cite{Kimura2016}.
By inequality (4) in \cite{Kimura2016}, and chebychev's inequality or Markov's inequality, we have 
\begin{align*}
b_n \max_{k} |t_{k}| = O_p(1).  
\end{align*}
Therefore, the typical term converges to $0$ in probability. 
The other terms are evaluated in the same way. 
We obtain the conclusion. 
%
%
\qed
\begin{proposition*} \label{prop1i}
$
\hat{\gamma}^{p q}_{n,i} = 
G_{n,i}^{p q} + o_p(b_n^{-1/2})$ 
under [B] or [B'], and [C] (i=1,2).  
\end{proposition*}
\proof
We just confirm the convergence of the easiest term. 
\begin{align*}
& b_n^{1+\rho} \sum_{k=2}^{b_n-2} |s_{k} \otimes s_{k} - t_{k} \otimes t_{k}| \\
&= b_n^{1+\rho} \sum_{k=2}^{b_n-2} 
\Big|(s_{k} - t_{k}) \otimes (s_{k} - t_{k}) 
+ (s_{k} - t_{k}) \otimes t_{k}
+ t_{k} \otimes (s_{k} - t_{k}) \Big| 
\end{align*}
By definition, 
\begin{align*}
s_k - t_k 
&= (\tilde{\mathbb{M}}^n_{k} - \tilde{\mathbb{M}}^n_{k-1})^{\otimes} 
+ (\tilde{\mathbb{M}}^n_{k} - \tilde{\mathbb{M}}^n_{k-1}) 
\otimes (\tilde{\bm{\chi}}^n_{k} - \tilde{\bm{\chi}}^n_{k-1}) \\
&\quad
+ (\tilde{\bm{\chi}}^n_{k} - \tilde{\bm{\chi}}^n_{k-1})
\otimes (\tilde{\mathbb{M}}^n_{k} - \tilde{\mathbb{M}}^n_{k-1}). 
\end{align*}
Therefore, it is enough to prove the following equations. 
\begin{align*}
b_n^{1+\rho} \sum_{k=2}^{b_n-2} 
(\tilde{\mathbb{M}}^n_{k} - \tilde{\mathbb{M}}^n_{k-1})^{\otimes 3}
&=o_p(1)\\
b_n^{1+\rho} \sum_{k=2}^{b_n-2} 
(\tilde{\mathbb{M}}^n_{k} - \tilde{\mathbb{M}}^n_{k-1})^{\otimes 2}
\otimes (\tilde{\bm{\chi}}^n_{k} - \tilde{\bm{\chi}}^n_{k-1})
&=o_p(1)\\
b_n^{1+\rho} \sum_{k=2}^{b_n-2} 
(\tilde{\mathbb{M}}^n_{k} - \tilde{\mathbb{M}}^n_{k-1})^{\otimes}
\otimes (\tilde{\bm{\chi}}^n_{k} - \tilde{\bm{\chi}}^n_{k-1})^{\otimes}
&=o_p(1)\\
b_n^{1+\rho} \sum_{k=2}^{b_n-2} 
(\tilde{\mathbb{M}}^n_{k} - \tilde{\mathbb{M}}^n_{k-1})
\otimes (\tilde{\bm{\chi}}^n_{k} - \tilde{\bm{\chi}}^n_{k-1})^{\otimes 2}
&=o_p(1), 
\end{align*}
where for vector $V$, 
$V^{\otimes} = V \otimes V$, 
$V^{\otimes 2} = V \otimes V \otimes V$, 
$V^{\otimes 3} = V \otimes V \otimes V \otimes V$. 
By assumption [C], it holds that
\begin{align*}
E\bigg[ \bigg| b_n^{1+\rho} \sum_{k=2}^{b_n-2} 
(\tilde{\mathbb{M}}^n_{k} - \tilde{\mathbb{M}}^n_{k-1})^{\otimes 3} \bigg| \bigg] 
&\leq 
E\bigg[ b_n^{1+\rho} \sum_{k=2}^{b_n-2} 
| \tilde{\mathbb{M}}^n_{k} - \tilde{\mathbb{M}}^n_{k-1} |^{4} \bigg] \\
&\leq 
16 E\bigg[ b_n^{1+\rho} \sum_{k=2}^{b_n-2} 
| \tilde{\mathbb{M}}^n_{k} |^{4} \bigg] \\
&= O(b_n^{4+\rho} a_n^{-2}) 
= o(1).
\end{align*}
Lemma 1 (b) and inequality (4) in \cite{Kimura2016} yields  
\begin{align*}
&\bigg| b_n^{1+\rho} \sum_{k=2}^{b_n-2} 
(\tilde{\mathbb{M}}^n_{k} - \tilde{\mathbb{M}}^n_{k-1})^{\otimes}
\otimes (\tilde{\bm{\chi}}^n_{k} - \tilde{\bm{\chi}}^n_{k-1})^{\otimes} \bigg| \\
&\leq 
\bigg( b_n^{\rho} \sum_{k=2}^{b_n-2} 
| \tilde{\mathbb{M}}^n_{k} - \tilde{\mathbb{M}}^n_{k-1} |^2 \bigg)
\bigg( b_n |\tilde{\bm{\chi}}^n_{k} - \tilde{\bm{\chi}}^n_{k-1} |^2 \bigg) 
= o_p(1) O_p(1) 
= o_p(1). 
\end{align*}
and
\begin{align*}
&\bigg| b_n^{1+\rho} \sum_{k=2}^{b_n-2} 
(\tilde{\mathbb{M}}^n_{k} - \tilde{\mathbb{M}}^n_{k-1})
\otimes (\tilde{\bm{\chi}}^n_{k} - \tilde{\bm{\chi}}^n_{k-1})^{\otimes \otimes} \bigg|  \\
&\leq 
\bigg| b_n^{\rho} \sum_{k=2}^{b_n-2} 
 (\tilde{\mathbb{M}}^n_{k} - \tilde{\mathbb{M}}^n_{k-1}) \otimes (\tilde{\bm{\chi}}^n_{k} - \tilde{\bm{\chi}}^n_{k-1}) \bigg|
\bigg| b_n (\tilde{\bm{\chi}}^n_{k} - \tilde{\bm{\chi}}^n_{k-1} )^{\otimes} \bigg| \\
&= o_p(1) O_p(1) 
= o_p(1). 
\end{align*}
From the above  inequalities, it is derived that 
\begin{align*}
&\bigg| b_n^{1+\rho} \sum_{k=2}^{b_n-2} 
(\tilde{\mathbb{M}}^n_{k} - \tilde{\mathbb{M}}^n_{k-1})^{\otimes 2}
\otimes (\tilde{\bm{\chi}}^n_{k} - \tilde{\bm{\chi}}^n_{k-1}) \bigg| \\
&\leq
\bigg| 
b_n^{1+\rho} \sum_{k=2}^{b_n-2} 
(\tilde{\mathbb{M}}^n_{k} - \tilde{\mathbb{M}}^n_{k-1})^{\otimes 3} 
\bigg|^{1/2} \\
&\qquad 
\times 
\bigg| b_n^{1+\rho} \sum_{k=2}^{b_n-2} 
(\tilde{\mathbb{M}}^n_{k} - \tilde{\mathbb{M}}^n_{k-1})^{\otimes}
\otimes (\tilde{\bm{\chi}}^n_{k} - \tilde{\bm{\chi}}^n_{k-1})^{\otimes} 
\bigg|^{1/2} \\
&= o_p(1) o_p(1) = o_p(1). 
\end{align*}
\qed

Now, we let 
\begin{align*}
A^\alpha_k &= 
X_{t_{k-2}}^{\alpha 1} \int_{I_{k-1}} H_k(s)  \ dw_s +X_{t_{k-2}}^{\alpha 1} \int_{I_{k}} K_k(s) \ dw_s \\
B^\alpha_k &= \frac{1}{h_{k}} \int_{I_{k}} \int_{t_{k-2}}^{t} (X_{s}^{\alpha 1}-X_{t_{k-2}}^{\alpha 1}) \ dw_s \ dt - \frac{1}{h_{k-1}} \int_{I_{k-1}}\int_{t_{k-2}}^{t}(X_{s}^{\alpha 1}-X_{t_{k-2}}^{\alpha 1}) \ dw_s  \ dt
\\
C^\alpha_k &= \frac{1}{h_{k}} \int_{I_{k}} \int_{t_{k-2}}^{t} X_{t_{k-2}}^{\alpha 0}ds \ dt - \frac{1}{h_{k-1}} \int_{I_{k-1}}\int_{t_{k-2}}^{t}X_{t_{k-2}}^{\alpha 0}ds  \ dt
\\
D^\alpha_k &= \frac{1}{h_{k}} \int_{I_{k}} \int_{t_{k-2}}^{t} (X_{s}^{\alpha 0}-X_{t_{k-2}}^{\alpha 0})ds \ dt - \frac{1}{h_{k-1}} \int_{I_{k-1}}\int_{t_{k-2}}^{t}(X_{s}^{\alpha 0}-X_{t_{k-2}}^{\alpha 0})ds  \ dt. 
\end{align*}
By Lemma 2 and inequality (4) in \cite{Kimura2016}, we have
\begin{align*}
G_{n,1}^{p q} 
&= \frac{9}{8} \Big\{ \sum_{k=2}^{b_n} 
(A^{\alpha_1}_{k} A^{\beta_1}_{k}) (A^{\alpha_2}_{k} A^{\beta_2}_{k}) \\
&\qquad\qquad\qquad
- \frac{1}{2} \sum_{k=2}^{b_n-2} 
\Big( (A^{\alpha_1}_{k} A^{\beta_1}_{k}) (A^{\alpha_2}_{k+2} A^{\beta_2}_{k+2}) \\
&\qquad\qquad\qquad\qquad\qquad
+ (A^{\alpha_1}_{k+2} A^{\beta_1}_{k+2}) (A^{\alpha_2}_{k} A^{\beta_2}_{k}) \Big)
\Big\} \Big(\frac{T}{b_n} \Big)^{-1} 
+ o_p(b_n^{-1/2}) \\
&= 
\mathcal{R}_{n,1}^{1,p q} + o_p(b_n^{-1/2}) 
\end{align*}
and
\begin{align*}
G_{n,2}^{p q} 
&= \frac{9}{8} \sum_{k=2}^{b_n-2} \frac{1}{2} \Big\{
(A^{\alpha_1}_{k+2} A^{\beta_1}_{k+2}) (A^{\alpha_2}_{k+2} A^{\beta_2}_{k+2}) \\
&\qquad\qquad\qquad\qquad\quad
- (A^{\alpha_1}_{k+2} A^{\beta_1}_{k+2}) (A^{\alpha_2}_{k} A^{\beta_2}_{k})  \\
&\qquad\qquad\qquad\qquad\quad
- (A^{\alpha_1}_{k} A^{\beta_1}_{k}) (A^{\alpha_2}_{k+2} A^{\beta_2}_{k+2})\\
&\qquad\qquad\qquad\qquad\quad
+ (A^{\alpha_1}_{k} A^{\beta_1}_{k}) (A^{\alpha_2}_{k} A^{\beta_2}_{k})
\Big\} \Big(\frac{T}{b_n} \Big)^{-1} + o_p(b_n^{-1/2}) \\
&= 
\mathcal{R}_{n,2}^{1,p q} + o_p(b_n^{-1/2}). 
\end{align*}
Here, we write 
\begin{align*}
q^{2,\alpha \beta}_k 
&= 
X^{\alpha1}_{t_{k-2}}\cdot X^{\beta1}_{t_{k-2}}
 \frac{h_{k-1}}{3} 
 + 
 X^{\alpha1}_{t_{k-2}}\cdot X^{\beta1}_{t_{k-2}}
 \frac{h_{k}}{3}
\end{align*}
and 
\begin{align*}
q_k^{3, \alpha \beta} 
&=
2 (\mathbb{X}_{t_{k-2}}^{1})^{\tilde{\otimes}(\alpha,\beta)}
\cdot
 \Big( \int_{I_{k-1}} \int_{t_{k-2}}^{t} H_k(s) \ dw_s 
 \otimes H_k(t) \ dw_t \\
&\qquad\qquad\qquad\qquad 
+ \int_{I_{k}} \int_{t_{k-1}}^{t} K_k(s) \ dw_s  \otimes K_k(t) \ dw_t \\
&\qquad\qquad\qquad\qquad
+ \int_{I_{k-1}} H_k(s) \ dw_s \otimes \int_{I_{k}}  K_k(s) \ dw_s \Big). 
\end{align*}
Then, we have
\begin{align*}
\mathcal{R}_{n,1}^{1,p q} 
&= \frac{9}{8} \Big\{ \sum_{k=2}^{b_n} 
(q_{k}^{2, p} + q_{k}^{3, p}) (q_{k}^{2, q} + q_{k}^{3, q}) \\
&\qquad\qquad\qquad
- \frac{1}{2} \sum_{k=2}^{b_n-2} 
\Big( (q_{k}^{2, p} + q_{k}^{3, p}) (q_{k+2}^{2, q} + q_{k+2}^{3, q}) \\
&\qquad\qquad\qquad\qquad\qquad
+ (q_{k+2}^{2, p} + q_{k+2}^{3, p}) (q_{k}^{2, q} + q_{k}^{3, q}) \Big)
\Big\} \Big(\frac{T}{b_n} \Big)^{-1} 
\end{align*}
and
\begin{align*}
\mathcal{R}_{n,2}^{1,p q} 
&= \frac{9}{8} \sum_{k=2}^{b_n-2} \frac{1}{2} \Big\{
(q_{k+2}^{2, p} + q_{k+2}^{3, p}) (q_{k+2}^{2, p} + q_{k+2}^{3, p}) \\
&\qquad\qquad\qquad\qquad
- (q_{k+2}^{2, p} + q_{k+2}^{3, p}) (q_{k}^{2, p} + q_{k}^{3, p}) \\
&\qquad\qquad\qquad\qquad
- (q_{k}^{2, p} + q_{k}^{3, p}) (q_{k+2}^{2, p} + q_{k+2}^{3, p}) \\
&\qquad\qquad\qquad\qquad
+ (q_{k}^{2, p} + q_{k}^{3, p}) (q_{k}^{2, p} + q_{k}^{3, p})
\Big\} \Big(\frac{T}{b_n} \Big)^{-1} 
\end{align*}
by It\^o's formula and simple calculus. 
We only see the typical terms. 
We have 
\begin{align} \label{2k2k}
\Big(\frac{T}{b_n} \Big)^{-1} \sum_{k=2}^{b_n} q_{k}^{2, p} q_{k}^{2, q} 
= \frac{4}{9} \sum_{k=2}^{b_n} (X_{t_{k-2}}^{\alpha_1 1} \cdot X_{t_{k-2}}^{\beta_1 1})
(X_{t_{k-2}}^{\alpha_2 1} \cdot X_{t_{k-2}}^{\beta_2 1}) \Big(\frac{T}{b_n} \Big), 
\end{align}
%
\begin{align} \label{2k3k}
\Big(\frac{T}{b_n} \Big)^{-1} \sum_{k=2}^{b_n} q_{k}^{2, p} q_{k}^{3, q} 
=\sum_{k=2}^{b_n} \frac{2}{3} (X_{t_{k-2}}^{\alpha_1 1} \cdot X_{t_{k-2}}^{\beta_1 1}) q_{k}^{3, q} 
= O_p(b_n^{-1/2}), 
\end{align}
%
and
\begin{align} \label{3k+d3k}
\Big(\frac{T}{b_n} \Big)^{-1} \sum_{k=2}^{b_n-2} q_{k+2}^{3, p} q_{k}^{3, q} 
= O_p(b_n^{-1/2}). 
\end{align}
By Lemma 4 in \cite{Kimura2016} and simple calculus, it holds that
\begin{align} \label{3k3k}
&\Big(\frac{T}{b_n} \Big)^{-1} \sum_{k=2}^{b_n} q_{k}^{3, p} q_{k}^{3, q} \notag\\
&= \Big(\frac{T}{b_n} \Big)^{-1} \sum_{k=2}^{b_n} 
\Big\{\sum_{i,j = 1}^{r}
(X^{\alpha_1 1}_{i, t_{k-2}} X^{\beta_1 1}_{j, t_{k-2}} + X^{\alpha_1 1}_{j, t_{k-2}} X^{\beta_1 1}_{i, t_{k-2}}) \notag\\
&\qquad\qquad\qquad\qquad\quad
\times 
\Big(\int_{I_{k-1}} \int_{t_{k-2}}^{t} H_k(s) \ dw^{i}_s H_k(t) \ dw^{j}_t \notag\\
&\qquad\qquad\qquad\qquad\qquad 
+ \int_{I_{k}} \int_{t_{k-1}}^{t} K_k(s) \ dw^{i}_s  K_k(t) \ dw^{j}_t \notag\\
&\qquad\qquad\qquad\qquad\qquad
+ \int_{I_{k-1}} H_k(s) \ dw^{i}_s \int_{I_{k}}  K_k(s) \ dw^{j}_s \Big) \Big\} \notag\\
&\qquad\qquad\qquad
\times 
\Big\{ \sum_{i,j = 1}^{r}
(X^{\alpha_2 1}_{i, t_{k-2}} X^{\beta_2 1}_{j, t_{k-2}} + X^{\alpha_2 1}_{j, t_{k-2}} X^{\beta_2 1}_{i, t_{k-2}}) \notag\\
&\qquad\qquad\qquad\qquad\quad
\times 
\Big(\int_{I_{k-1}} \int_{t_{k-2}}^{t} H_k(s) \ dw^{i}_s H_k(t) \ dw^{j}_t \notag\\
&\qquad\qquad\qquad\qquad\qquad 
+ \int_{I_{k}} \int_{t_{k-1}}^{t} K_k(s) \ dw^{i}_s  K_k(t) \ dw^{j}_t \notag\\
&\qquad\qquad\qquad\qquad\qquad 
+ \int_{I_{k-1}} H_k(s) \ dw^{i}_s \int_{I_{k}}  K_k(s) \ dw^{j}_s \Big) \Big\} \notag\\
&= 
\Big(\frac{T}{b_n} \Big)^{-1} \sum_{k=2}^{b_n} 
\Big\{\sum_{i,j = 1}^{r}
(X^{\alpha_1 1}_{i, t_{k-2}} X^{\beta_1 1}_{j, t_{k-2}} + X^{\alpha_1 1}_{j, t_{k-2}} X^{\beta_1 1}_{i, t_{k-2}}) \notag\\
&\qquad\qquad\qquad\qquad\quad 
\times
(X^{\alpha_2 1}_{i, t_{k-2}} X^{\beta_2 1}_{j, t_{k-2}} + X^{\alpha_2 1}_{j, t_{k-2}} X^{\beta_2 1}_{i, t_{k-2}}) \Big\} \notag\\
&\qquad\qquad\qquad 
\times 
\Big\{ \int_0^T \int_0^t (\tilde{H}_k(s))^2 \ ds (\tilde{H}_k(t))^2 \ dt \notag\\
&\qquad\qquad\qquad\qquad 
+ \int_0^T \int_0^t (\tilde{K}_k(s))^2 \ ds (\tilde{K}_k(t))^2 \ dt \notag\\
&\qquad\qquad\qquad\qquad 
+ \int_0^T (\tilde{H}_k(s))^2 \ ds \int_0^T (\tilde{H}_k(s))^2 \ ds \Big\} 
+ O_p(b_n^{-1/2}) \notag\\
&= \frac{4}{9} \sum_{k=2}^{b_n} \Big\{
(X^{\alpha_1 1}_{t_{k-2}} \cdot X^{\alpha_2 1}_{t_{k-2}})(X^{\beta_1 1}_{t_{k-2}} \cdot X^{\beta_2 1}_{t_{k-2}}) \notag\\
&\qquad\qquad
+ (X^{\alpha_1 1}_{t_{k-2}} \cdot X^{\beta_2 1}_{t_{k-2}})(X^{\beta_1 1}_{t_{k-2}} \cdot X^{\alpha_2 1}_{t_{k-2}}) \Big\} 
\Big(\frac{T}{b_n}\Big) 
+ O_p(b_n^{-1/2}), 
\end{align}
where $\tilde{H}_k(s) = 1_{I_{k-1}}(s) H_k(s)$ and $\tilde{K}_k(s) = 1_{I_{k}}(s) K_k(s)$. 
%
Therefore, we have 
\begin{align*}
\mathcal{R}_{n,1}^{1,p q} 
&= \frac{9}{8} \bigg\{ \frac{4}{9} \sum_{k=2}^{b_n} (X_{t_{k-2}}^{\alpha_1 1} \cdot X_{t_{k-2}}^{\beta_1 1})
(X_{t_{k-2}}^{\alpha_2 1} \cdot X_{t_{k-2}}^{\beta_2 1}) \Big(\frac{T}{b_n}\Big) \\
&\qquad
+ \frac{4}{9} \sum_{k=2}^{b_n} \Big\{
(X^{\alpha_1 1}_{t_{k-2}} \cdot X^{\alpha_2 1}_{t_{k-2}})(X^{\beta_1 1}_{t_{k-2}} \cdot X^{\beta_2 1}_{t_{k-2}}) \\
&\qquad\qquad\qquad
+ (X^{\alpha_1 1}_{t_{k-2}} \cdot X^{\beta_2 1}_{t_{k-2}})(X^{\beta_1 1}_{t_{k-2}} \cdot X^{\alpha_2 1}_{t_{k-2}}) \Big\} 
(\frac{T}{b_n}\Big) \\
&\qquad
- \frac{1}{2} \Big\{ \frac{4}{9} \sum_{k=2}^{b_n-2} (X_{t_{k-2}}^{\alpha_1 1} \cdot X_{t_{k-2}}^{\beta_1 1})
(X_{t_{k}}^{\alpha_2 1} \cdot X_{t_{k}}^{\beta_2 1}) \Big(\frac{T}{b_n}\Big) \\
&\qquad\qquad
+ \frac{4}{9} \sum_{k=2}^{b_n-2} (X_{t_{k}}^{\alpha_1 1} \cdot X_{t_{k}}^{\beta_1 1})
(X_{t_{k-2}}^{\alpha_2 1} \cdot X_{t_{k-2}}^{\beta_2 1}) \Big(\frac{T}{b_n}\Big) \Big\} \bigg\} \\
&\qquad
+ O_p(b_n^{-1/2}) \\
&= \gamma^{p q} + o_p(b_n^{-1/2}) + O_p(b_n^{-1/2}) 
\end{align*}
and
\begin{align*}
\mathcal{R}_{n,2}^{1,p q} 
&= \frac{9}{16} 
\bigg\{
\frac{4}{9} \sum_{k=2}^{b_n-2} (X_{t_{k}}^{\alpha_1 1} \cdot X_{t_{k}}^{\beta_1 1})
(X_{t_{k}}^{\alpha_2 1} \cdot X_{t_{k}}^{\beta_2 1}) \Big(\frac{T}{b_n}\Big) \\
&\qquad\quad
+ \frac{4}{9} \sum_{k=2}^{b_n-2} \Big\{
(X^{\alpha_1 1}_{t_{k-2}} \cdot X^{\alpha_2 1}_{t_{k-2}})(X^{\beta_1 1}_{t_{k-2}} \cdot X^{\beta_2 1}_{t_{k-2}}) \\
&\qquad\qquad\qquad\quad
+ (X^{\alpha_1 1}_{t_{k-2}} \cdot X^{\beta_2 1}_{t_{k-2}})(X^{\beta_1 1}_{t_{k-2}} \cdot X^{\alpha_2 1}_{t_{k-2}}) \Big\} 
\Big(\frac{T}{b_n}\Big) \\
&\qquad\quad
- \frac{4}{9} \sum_{k=2}^{b_n-2} (X_{t_{k}}^{\alpha_1 1} \cdot X_{t_{k}}^{\beta_1 1})
(X_{t_{k-2}}^{\alpha_2 1} \cdot X_{t_{k-2}}^{\beta_2 1}) \Big(\frac{T}{b_n}\Big) \\
&\qquad\quad 
- \frac{4}{9} \sum_{k=2}^{b_n-2} (X_{t_{k-2}}^{\alpha_1 1} \cdot X_{t_{k-2}}^{\beta_1 1})
(X_{t_{k}}^{\alpha_2 1} \cdot X_{t_{k}}^{\beta_2 1}) \Big(\frac{T}{b_n}\Big) \\
&\qquad\quad 
+ \frac{4}{9} \sum_{k=2}^{b_n-2} (X_{t_{k-2}}^{\alpha_1 1} \cdot X_{t_{k-2}}^{\beta_1 1})
(X_{t_{k-2}}^{\alpha_2 1} \cdot X_{t_{k-2}}^{\beta_2 1}) \Big(\frac{T}{b_n}\Big) \\
&\qquad\quad 
+ \frac{4}{9} \sum_{k=2}^{b_n-2} \Big\{
(X^{\alpha_1 1}_{t_{k-2}} \cdot X^{\alpha_2 1}_{t_{k-2}})(X^{\beta_1 1}_{t_{k-2}} \cdot X^{\beta_2 1}_{t_{k-2}})  \\
&\qquad\qquad\qquad\quad 
+ (X^{\alpha_1 1}_{t_{k-2}} \cdot X^{\beta_2 1}_{t_{k-2}})(X^{\beta_1 1}_{t_{k-2}} \cdot X^{\alpha_2 1}_{t_{k-2}}) \Big\} 
\Big(\frac{T}{b_n}\Big) \bigg\} \\
&\qquad
+ O_p(b_n^{-1/2}). \\
&= \gamma^{p q} + o_p(b_n^{-1/2}) + O_p(b_n^{-1/2}). 
\end{align*}
The last $o_p(b_n^{-1/2})$s are derived from the proof of Lemma 3 in \cite{Kimura2016} and simple calculus. 
\qed

\subsection*{Proof of Theorem \ref{kernel}}
%
In the similar way to the proof of Therem \ref{nonkernel}, let  
\begin{align*}
G_{n,h}^{pq} 
&= \frac{9}{8} \sum_{k=2}^{b_n} \Big\{ 
\partial_{h} \{ \tilde{\chi}^{n, \alpha_1} , \tilde{\chi}^{n, \alpha_2} \}_{k} \
\partial_{h} \{ \tilde{\chi}^{n, \beta_1} , \tilde{\chi}^{n, \beta_2} \}_{k} \\
&\qquad\qquad
+ \partial_{h} \{ \tilde{\chi}^{n, \alpha_1} , \tilde{\chi}^{n, \beta_2} \}_{k} \
\partial_{h} \{ \tilde{\chi}^{n, \beta_1} , \tilde{\chi}^{n, \alpha_2} \}_{k} \Big\} \Big(\frac{T}{b_n} \Big). 
\end{align*}
%
%
%
%
%
\begin{lemma*} \label{lemma1h}
(a)
$
\hat{\gamma}^{p q}_{n,h} = 
G_{n,h}^{p q} + o_p(1)$ 
under any one of [B], [B$'$] and [B$^{\sharp^{\sharp}}$]. 

(b)
$
\hat{\gamma}^{p q}_{n,h} = 
G_{n,h}^{p q} + o_p(b_n^{-1/2})$ 
under [B$^{\sharp^{\sharp}}$]. 
\end{lemma*}

\proof
%
%
We write
\begin{align*}
\check{s}_{k, h}
:= \sum_{l=(k - n(h) +1) \vee 2}^{k} s_{l}
\ h^{-1}, 
\qquad
\check{t}_{k, h}
:= \sum_{l=(k - n(h) +1) \vee 2}^{k} t_{l}
\ h^{-1},  
\end{align*}
$\hat{\gamma}_{n,h} = (\hat{\gamma}^{p q}_{n,h})_{p, q = (1,2),(2,1),(1,1),(2,2)}$, and 
${G}_{n,h} = ({G}^{p q}_{n,h})_{p, q = (1,2),(2,1),(1,1),(2,2)}$. 
It holds that 
\begin{align*}
&b_n^{\rho} \Big(\frac{9}{8} \Big)^{-1} 
| \hat{\gamma}_{n,h} - G_{n,h} | \\
&= b_n^{\rho} 2 
\bigg| \sum_{k=2}^{b_n} (\check{s}_{k, h} \otimes \check{s}_{k, h} - \check{t}_{k, h} \otimes \check{t}_{k, h}) \bigg| \Big(\frac{T}{b_n} \Big) \\
&= b_n^{\rho} 2 
\bigg| \sum_{k=2}^{b_n} (\check{s}_{k, h} - \check{t}_{k, h}) \otimes (\check{s}_{k, h} - \check{t}_{k, h}) \\
&\qquad\qquad\qquad\qquad
+ (\check{s}_{k, h} - \check{t}_{k, h}) \otimes \check{t}_{k, h} 
+ \check{t}_{k, h} \otimes (\check{s}_{k, h} - \check{t}_{k, h}) \bigg| \Big(\frac{T}{b_n} \Big) \\
&\leq b_n^{\rho} 2 
\sum_{k=2}^{b_n} 
\Big\{ \Big| (\check{s}_{k, h} - \check{t}_{k, h}) \otimes (\check{s}_{k, h} - \check{t}_{k, h}) \Big| \\
&\qquad\qquad\qquad\qquad
+ \Big| (\check{s}_{k, h} - \check{t}_{k, h}) \otimes \check{t}_{k, h} \Big| 
+ \Big| \check{t}_{k, h} \otimes (\check{s}_{k, h} - \check{t}_{k, h}) \Big| \Big\} 
\Big(\frac{T}{b_n} \Big) \\
&= b_n^{\rho} 2 
\sum_{k=2}^{b_n} | \check{s}_{k, h} - \check{t}_{k, h} |^2 \Big(\frac{T}{b_n} \Big) 
+ b_n^{\rho} 4 
\sum_{k=2}^{b_n} | \check{t}_{k, h} | | \check{s}_{k, h} - \check{t}_{k, h} | \Big(\frac{T}{b_n}\Big) \\
&\leq 2 \bigg\{ b_n^{\rho/2} \Big(\frac{T}{b_n} \Big)^{-1/2} \sum_{k=2}^{b_n} | \check{s}_{k, h} - \check{t}_{k, h} | \Big(\frac{T}{b_n} \Big) \bigg\}^2 \\
&\quad
+ 4 \Big\{ \max_k |\check{t}_{k, h} | \Big\}
\bigg\{ b_n^{\rho} \sum_{k=2}^{b_n} | \check{s}_{k, h} - \check{t}_{k, h} | \Big(\frac{T}{b_n} \Big) \bigg\}. 
\end{align*}
We obtain the conclusion from the inequalities
\begin{align*}
\sum_{k=2}^{b_n} | \check{s}_{k, h} - \check{t}_{k, h} | \Big(\frac{T}{b_n} \Big) 
\leq \sum_{k=2}^{b_n} \sum_{l=(k - n(h) +1) \vee 2}^{k} | s_{l} - t_{l} | \ h^{-1} \Big(\frac{T}{b_n} \Big)  
\leq \sum_{k=2}^{b_n} | s_{k} - t_{k} |, 
\end{align*}
\begin{align*}
\max_{k} | \check{t}_{k, h} | \leq \max_{k} | t_{k} |,
\end{align*} 
and the proof of Lemma \ref{lemma1i}. 
\qed
%
%
%
\begin{proposition*} \label{prop1h}
$
\hat{\gamma}^{p q}_{n,h} = 
G_{n,h}^{p q} + o_p(b_n^{-1/2})$ 
under [B] or [B'], and [C].  
\end{proposition*}
\proof
The conclusion is derived from the proof of Lemma \ref{lemma1h}, that of Proposition \ref{prop1i}, and the inequality
\begin{align*}
b_n^{\rho} 
\sum_{k=2}^{b_n} | \check{s}_{k, h} - \check{t}_{k, h} |^2 \Big(\frac{T}{b_n} \Big)
\leq
b_n^{\rho} \sum_{k=2}^{b_n} | {s}_{k} - {t}_{k} |^2 \Big(\frac{T}{b_n} \Big)^{-1}.
\end{align*} 
\qed
\begin{lemma*} \label{lemma2h}
$b_n^{\rho} \sum_{l=(k - n(h) +1) \vee 2}^{k} F_{l} G_{l} h^{-1} \to^p 0$
for all pairs of  \\
$(F_{l}^, G_{l}) \in \{ A_{k}^{\alpha}, B_{k}^{\alpha}, C_{k}^{\alpha}, D_{k}^{\alpha} ; \alpha=1, 2 \}^2 \backslash \{(A_{l}^{\alpha}, A_{l}^{\beta}) ; \alpha, \beta =1, 2 \}$.
\end{lemma*}

The proof of this lemma is very similar to proof of Lemma 2 in \cite{Kimura2016}. 

\proof
Let 
$\mathcal{Q}_{n, h} = b_n^{\rho} \sum_{l=(k - n(h) +1) \vee 2}^{k} F_{l} G_{l} h^{-1}$. 
%
%
For $F_{l} = A_{l}^{\alpha}$ and $G_{l} = C_{l}^{\beta}$, 
$E[(|\mathcal{Q}_{n, h}|)^2] = O(b_n^{2 \rho} n(h) h^{-2} b_n^{-3}) = O(h^{-1} b_n^{2\rho -2})=o(1)$. 
For $F_{l} = A_{l}^{\alpha}$ and $G_{l} = D_{l}^{\beta}$, 
$E[|\mathcal{Q}_{n, h}|] = O( b_n^{\rho} n(h) h^{-1} b_n^{-3/2}) o(1) = o(b_n^{\rho -1/2})=o(1)$. 
%
%
%
%
%
For $F_{l} = A_{l}^{\alpha}$ and $G_{l} = B_{l}^{\beta}$, 
$E[( |\mathcal{Q}_{n, h}|)^2] = O(b_n^{2 \rho} n(h) h^{-2} b_n^{-3}) = O(h^{-1} b_n^{2 \rho -2})=o(1)$. 

\qed
%
%
%
%

Now, let
\begin{align*}
Q_{k, h}^{1, \alpha \beta} = \sum_{l=(k - n(h) +1) \vee 2}^{k} A_{l}^{\alpha} A_{l}^{\beta} h^{-1}. 
\end{align*}
By Lemma \ref{lemma2h} 
and (4) 
(esp. $Q_{k, h}^{1, \alpha \beta} = O_p(1)$), 
it holds that 
\begin{align*}
G_{n,h}^{p q} 
&= 
\frac{9}{8} \sum_{k=2}^{b_n} (Q_{k, h}^{1, \alpha_1 \alpha_2} Q_{k, h}^{1, \beta_1 \beta_2} 
+ Q_{k, h}^{1, \alpha_1 \beta_2} Q_{k, h}^{1, \beta_1 \alpha_2}) \Big(\frac{T}{b_n} \Big) + o_p(b_n^{-1/2}) \\
&=: 
\mathcal{R}_{n,h}^{1, p q} +  o_p(b_n^{-1/2}). 
\end{align*}
%

By recombination of the multiple summation, 
\begin{align*}
\mathcal{R}_{n,h}^{1, p q}
&= \frac{9}{8}  \sum_{l=2}^{b_n - n(h)} \Big\{
n(h) 
(q_{l}^{1, \alpha_1 \alpha_2} q_{l}^{1, \beta_1 \beta_2} 
+ q_{l}^{1, \alpha_1 \beta_2} q_{l}^{1, \beta_1 \alpha_2}) \\
&\qquad\qquad\quad 
+ (n(h) - 1) 
(q_{l}^{1, \alpha_1 \alpha_2} q_{l+1}^{1, \beta_1 \beta_2} 
+ q_{l}^{1, \alpha_1 \beta_2} q_{l+1}^{1, \beta_1 \alpha_2} \\
&\qquad\qquad\qquad\qquad\qquad 
+ q_{l+1}^{1, \alpha_1 \alpha_2} q_{l}^{1, \beta_1 \beta_2} 
+ q_{l+1}^{1, \alpha_1 \beta_2} q_{l}^{1, \beta_1 \alpha_2}) \\
&\qquad\qquad\quad 
+ \sum_{k=2}^{n(h)-1} (n(h) - k) 
(q_{l}^{1, \alpha_1 \alpha_2} q_{l+k}^{1, \beta_1 \beta_2} 
+ q_{l}^{1, \alpha_1 \beta_2} q_{l+k}^{1, \beta_1 \alpha_2} \\
&\qquad\qquad\qquad\qquad\qquad\qquad\quad 
+ q_{l+k}^{1, \alpha_1 \alpha_2} q_{l}^{1, \beta_1 \beta_2} 
+ q_{l+k}^{1, \alpha_1 \beta_2} q_{l}^{1, \beta_1 \alpha_2})\Big\} h_n^{-2}  \Big(\frac{T}{b_n} \Big) \\
&\quad
+ \frac{9}{8}  \sum_{m=b_n - n(h) + 1}^{b_n} \Big\{
\Big( \sum_{l = b_n - n(h) + 1}^{m} q_{l}^{1, \alpha_1 \alpha_2} \Big)
\Big( \sum_{l = b_n - n(h) + 1}^{m} q_{l}^{1, \beta_1 \beta_2} \Big) \\
&\qquad\qquad\qquad\qquad
+ 
\Big( \sum_{l = b_n - n(h) + 1}^{m} q_{l}^{1, \alpha_1 \beta_2} \Big)
\Big( \sum_{l = b_n - n(h) + 1}^{m} q_{l}^{1, \beta_1 \alpha_2} \Big) \Big\} h_n^{-2}  \Big(\frac{T}{b_n} \Big) \\
&=: \mathcal{B}_{n,h}^{1, p q} + \mathcal{T}_{n,h}^{1, p q}, 
\end{align*}
where $q_{l}^{1, \alpha \beta} = A_{l}^{\alpha} A_{l}^{\beta} = q_{l}^{2, \alpha \beta} + q_{l}^{3, \alpha \beta}$. 

By (\ref{2k3k}), (\ref{3k+d3k}), (\ref{3k3k}) and similar evaluation, we have
$\mathcal{B}_{n,h}^{1, p q} = \mathcal{B}_{n,h}^{2, p q} + O_p(b_n^{-1/2})$, 
where $\mathcal{B}_{n,h}^{2, p q}$ is defined in the same way of  $\mathcal{B}_{n,h}^{1, p q}$ 
with  $q_{l}^{2, \alpha \beta}$ instead of $q_{l}^{1, \alpha \beta}$.
It is easy to see that 
$\mathcal{T}_{n,h}^{1, p q} = O_p(h_n)$. 

By simple calculus, 
we have
\begin{align*}
\mathcal{B}_{n,h}^{2, p q} 
&= 
\frac{9}{8} \sum_{l=2}^{b_n - n(h)} \sum_{k=1}^{n(h) -1} (n(h) - k) \\
&\qquad 
\times \Big\{
(q_{l}^{2, \alpha_1 \alpha_2} (q_{l+k}^{2, \beta_1 \beta_2} - q_{l}^{2, \beta_1 \beta_2})
+ q_{l}^{2, \alpha_1 \beta_2} (q_{l+k}^{2, \beta_1 \alpha_2} - q_{l}^{2, \beta_1 \alpha_2}) \\
&\qquad\quad 
+ (q_{l+k}^{2, \alpha_1 \alpha_2} - q_{l}^{2, \alpha_1 \alpha_2}) q_{l}^{2, \beta_1 \beta_2} 
+ (q_{l+k}^{2, \alpha_1 \beta_2} - q_{l}^{2, \alpha_1 \beta_2}) q_{l}^{2, \beta_1 \alpha_2}) 
\Big\} h_n^{-2}  \Big(\frac{T}{b_n} \Big) \\
&\quad
+ \frac{9}{8} \sum_{l=2}^{b_n} n(h)^2 
(q_{l}^{2, \alpha_1 \alpha_2} q_{l}^{2, \beta_1 \beta_2} 
+ q_{l}^{2, \alpha_1 \beta_2} q_{l}^{2, \beta_1 \alpha_2}) 
h_n^{-2}  \Big(\frac{T}{b_n} \Big) \\
&\quad
- \frac{9}{8} \sum_{l=b_n - n(h)+1}^{b_n} n(h)^2 
(q_{l}^{2, \alpha_1 \alpha_2} q_{l}^{2, \beta_1 \beta_2} 
+ q_{l}^{2, \alpha_1 \beta_2} q_{l}^{2, \beta_1 \alpha_2}) 
h_n^{-2}  \Big(\frac{T}{b_n} \Big) \\
&=: \mathcal{P}_{n,h}^{2, p q} 
+ \mathcal{C}_{n,h}^{2, p q}  
- \mathcal{T}_{n,h}^{2, p q}. 
\end{align*}
From the proof of Lemma 3 in \cite{Kimura2016} and simple calculus, it is derived that
$\mathcal{C}_{n,h}^{2, p q} = \gamma^{p q} + o_p(b_n^{-1/2})$. 
It is easy to see that 
$\mathcal{T}_{n,h}^{2, p q} = o_p(h_n)$. 
%
%

Furthermore, 
we have
\begin{align*}
\mathcal{P}_{n,h}^{2, p q} 
&= 
\frac{9}{8} \sum_{l=2}^{b_n - n(h)} \sum_{k=1}^{n(h) -1} (n(h) - k) \\
&\qquad 
\times \Big\{
(q_{l}^{2, \alpha_1 \alpha_2} \sum_{m=1}^{k} (q_{l+m}^{2, \beta_1 \beta_2} - q_{l+m-1}^{2, \beta_1 \beta_2})
+ q_{l}^{2, \alpha_1 \beta_2} \sum_{m=1}^{k} (q_{l+m}^{2, \beta_1 \alpha_2} - q_{l+m-1}^{2, \beta_1 \alpha_2}) \\
&\qquad\quad 
+ \sum_{m=1}^{k} (q_{l+m}^{2, \alpha_1 \alpha_2} - q_{l+m-1}^{2, \alpha_1 \alpha_2}) q_{l}^{2, \beta_1 \beta_2} 
+ \sum_{m=1}^{k} (q_{l+m}^{2, \alpha_1 \beta_2} - q_{l+m-1}^{2, \alpha_1 \beta_2}) q_{l}^{2, \beta_1 \alpha_2}) 
\Big\} h_n^{-2}  \Big(\frac{T}{b_n} \Big) \\
&= 
\frac{9}{8} \sum_{l=2}^{b_n - n(h)} \sum_{k=1}^{n(h) -1} \frac{(n(h) - k)(n(h) - k+1)}{2} \\
&\qquad 
\times \Big\{
(q_{l}^{2, \alpha_1 \alpha_2} (q_{l+k}^{2, \beta_1 \beta_2} - q_{l+k-1}^{2, \beta_1 \beta_2})
+ q_{l}^{2, \alpha_1 \beta_2} (q_{l+k}^{2, \beta_1 \alpha_2} - q_{l+k-1}^{2, \beta_1 \alpha_2}) \\
&\qquad\quad 
+ (q_{l+k}^{2, \alpha_1 \alpha_2} - q_{l+k-1}^{2, \alpha_1 \alpha_2}) q_{l}^{2, \beta_1 \beta_2} 
+ (q_{l+k}^{2, \alpha_1 \beta_2} - q_{l+k-1}^{2, \alpha_1 \beta_2}) q_{l}^{2, \beta_1 \alpha_2}) 
\Big\} h_n^{-2}  \Big(\frac{T}{b_n} \Big) \\
&= 
\frac{9}{8} \sum_{l=2}^{n(h)} \sum_{k=1}^{n(h) - l} \frac{(n(h) - k)(n(h) - k+1)}{2} \\
&\qquad 
\times \Big\{
(q_{l}^{2, \alpha_1 \alpha_2} (q_{l+k}^{2, \beta_1 \beta_2} - q_{l+k-1}^{2, \beta_1 \beta_2})
+ q_{l}^{2, \alpha_1 \beta_2} (q_{l+k}^{2, \beta_1 \alpha_2} - q_{l+k-1}^{2, \beta_1 \alpha_2}) \\
&\qquad\quad 
+ (q_{l+k}^{2, \alpha_1 \alpha_2} - q_{l+k-1}^{2, \alpha_1 \alpha_2}) q_{l}^{2, \beta_1 \beta_2} 
+ (q_{l+k}^{2, \alpha_1 \beta_2} - q_{l+k-1}^{2, \alpha_1 \beta_2}) q_{l}^{2, \beta_1 \alpha_2}) 
\Big\} h_n^{-2}  \Big(\frac{T}{b_n} \Big) \\
&\quad
+
\frac{9}{8} \sum_{l'=n(h)+1}^{b_n} \sum_{k=1}^{n(h) -1} \frac{(n(h) - k)(n(h) - k+1)}{2} \\
&\qquad 
\times \Big\{
(q_{l'-k}^{2, \alpha_1 \alpha_2} (q_{l'}^{2, \beta_1 \beta_2} - q_{l'-1}^{2, \beta_1 \beta_2})
+ q_{l'-k}^{2, \alpha_1 \beta_2} (q_{l'}^{2, \beta_1 \alpha_2} - q_{l'-1}^{2, \beta_1 \alpha_2}) \\
&\qquad\quad 
+ (q_{l'}^{2, \alpha_1 \alpha_2} - q_{l'-1}^{2, \alpha_1 \alpha_2}) q_{l'-k}^{2, \beta_1 \beta_2} 
+ (q_{l'}^{2, \alpha_1 \beta_2} - q_{l'-1}^{2, \alpha_1 \beta_2}) q_{l'-k}^{2, \beta_1 \alpha_2}) 
\Big\} h_n^{-2}  \Big(\frac{T}{b_n} \Big) \\
&\quad
-
\frac{9}{8} \sum_{l=b_n - n(h) + 1}^{b_n - 1} \sum_{k=1}^{b_n - l} \frac{(n(h) - k)(n(h) - k+1)}{2} \\
&\qquad 
\times \Big\{
(q_{l}^{2, \alpha_1 \alpha_2} (q_{l+k}^{2, \beta_1 \beta_2} - q_{l+k-1}^{2, \beta_1 \beta_2})
+ q_{l}^{2, \alpha_1 \beta_2} (q_{l+k}^{2, \beta_1 \alpha_2} - q_{l+k-1}^{2, \beta_1 \alpha_2}) \\
&\qquad\quad 
+ (q_{l+k}^{2, \alpha_1 \alpha_2} - q_{l+k-1}^{2, \alpha_1 \alpha_2}) q_{l}^{2, \beta_1 \beta_2} 
+ (q_{l+k}^{2, \alpha_1 \beta_2} - q_{l+k-1}^{2, \alpha_1 \beta_2}) q_{l}^{2, \beta_1 \alpha_2}) 
\Big\} h_n^{-2}  \Big(\frac{T}{b_n} \Big) \\
&=: 
\mathcal{H}_{n,h}^{3, p q}  
+ \mathcal{P}_{n,h}^{3, p q} 
- \mathcal{T}_{n,h}^{3, p q}. 
\end{align*}
$\mathcal{H}_{n,h}^{3, p q} = O_p(h_n)$ 
and 
$\mathcal{T}_{n,h}^{3, p q} = O_p(h_n)$
are derived from the following calculus. 
\begin{align*}
&\frac{(n(h) - k)(n(h) - k+1)}{2} (q_{l+k}^{2, \beta_1 \beta_2} - q_{l+k-1}^{2, \beta_1 \beta_2})\\
&= \frac{(n(h) - k)(n(h) - (k-1))}{2} q_{l+k}^{2, \beta_1 \beta_2}
- \frac{(n(h) - (k-1))(n(h) - (k-2))}{2} q_{l+(k-1)}^{2, \beta_1 \beta_2} \\
&\quad+ (n(h) - (k-1))q_{l+(k-1)}^{2, \beta_1 \beta_2}. 
\end{align*} 

Here, we set 
\begin{align*}
\mathcal{P}_{n,h}^{3, p q} 
= \sum_{l'=n(h)+1}^{b_n} p_{n,h, l'}^{3, p q}. 
\end{align*}
and we have 
\begin{align*}
\sum_{l'=n(h)+1}^{b_n} E\Big[p_{n,h, l'}^{3, p q} \Big| \mathcal{F}_{l'-3}^{n}\Big] &= O_p(h_n), \\ 
\sum_{l'=n(h)+1}^{b_n} E\Big[ \big|p_{n,h, l'}^{3, p q} \big|^2 \Big| \mathcal{F}_{l'-3}^{n}\Big] &= O_p(h_n^2) 
\end{align*}
and $\{p_{n,h, l'}^{3, p q}\}_{l'=2}^{b_n}$ is $\{\mathcal{F}_{l'-2}^{n}\}_{l'=2}^{b_n}$-adapted.  
Therefore, by Lenglart's inequality and its application, it holds that 
$\mathcal{P}_{n,h}^{3, p q} = O_p(h_n)$. 

Here, we proved that $G^{pq}_{n,h} = \gamma^{p q} + o_p(b_n^{-1/2}) +  O_p(b_n^{-1/2}) + O_p(h_n)$. 
We reached the conclusion.
\qed

\section{Simulation} \label{simulation}
We consider the following example.  
$\mathbb{X}$ is an $\mathbb{R}_{+}^2$-valued It\^o process satisfying
\begin{eqnarray*}
X^{1}_t 
&=& X^{1}_{0} + \int_{0}^{t} X^{1}_{s} \mu_{1} \ ds + \int_{0}^{t} X^{1}_{s} \sigma_{1} \ dw^{1}_{s} \\
X^{2}_t 
&=& X^{2}_{0} + \int_{0}^{t} X^{2}_{s} \mu_{2} \ ds + \int_{0}^{t} X^{2}_{s} \rho \sigma_{2} \ dw^{1}_{s} 
+ \int_{0}^{t} X^{2}_{s} \sqrt{1 - \rho^2} \sigma_{2} \ dw^{2}_{s},  
\end{eqnarray*}
where $w^1$ and $w^2$ are independent Wiener processes, and $t \in [0, T]$. 
$\mathbb{Y}$ is a counting process with intensity process $a_n \mathbb{X}$. 
In this setting, [C] is satisfied. 
 
In various $a_n$, $b_n$ settings, 
We examine the asymptotic variance estimators 
$\hat{\Xi}_{n,*}$, for $* = 1, 2, w, m, n$, where
$\hat{\Xi}_{n,w} =  \hat{\Xi}_{n,h}$, $h = h_n = T b_n^{-0.25}$; 
$\hat{\Xi}_{n,m} =  \hat{\Xi}_{n,h}$, $h = h_n = T b_n^{-0.5}$; 
$\hat{\Xi}_{n,n} =  \hat{\Xi}_{n,h}$, $h = h_n = T b_n^{-0.75}$.  
To reflect the condition that $\lim_{n \to \infty} b_n^r / a_n = 0$ ($r = 2, 2.5, 3$), 
we see the case that $a_n = b_n^r$ ($r = 2, 2.5, 3, 3.5$). 
We calculate these with 1000 paths and compare these by mean squared error (MSE). 

Now, we fix the parameters 
$\mu_1 = 0.2, 
\mu_2 = 0.3, 
\sigma_1 = 0.2, 
\sigma_2 = 0.3, 
\rho = 0.7, 
X^1_0 = 1, 
X^2_0 = 2, 
T = 1$. 

\begin{table}[tbp]
\centering
\caption{MSE of asymptotic variance estimators; $a_n = b_n^2$}
\begin{tabular}{rrrrrrrr}
  \hline
$*$ & $b_n = 2^4$ & $b_n = 2^5$ & $b_n = 2^6$ & $b_n = 2^7$ & $b_n = 2^8$ & $b_n = 2^9$ & $b_n = 2^{10}$ \\ 
  \hline
1 & 0.6514 & 0.6775 & 0.6931 & 0.6835 & 0.6892 & 0.6923 & 0.6908 \\ 
  2 & 0.4758 & 0.5838 & 0.6427 & 0.6537 & 0.6749 & 0.6845 & 0.6865 \\ 
  w & 0.2102 & 0.2279 & 0.2783 & 0.3554 & 0.3980 & 0.4213 & 0.4667 \\ 
  m & 0.6063 & 0.3973 & 0.7440 & 0.6470 & 0.7390 & 0.6298 & 0.7188 \\ 
  n & 1.3349 & 0.7355 & 0.3196 & 0.8940 & 1.3771 & 0.5679 & 0.6630 \\ 
   \hline
\end{tabular} \label{MSE2.0}
\end{table}

\begin{table}[tbp]
\centering
\caption{MSE of asymptotic variance estimators; $a_n = b_n^{2.5}$}
\begin{tabular}{rrrrrrrr}
  \hline
$*$ & $b_n = 2^4$ & $b_n = 2^5$ & $b_n = 2^6$ & $b_n = 2^7$ & $b_n = 2^8$ & $b_n = 2^9$ & $b_n = 2^{10}$ \\ 
  \hline
1 & 0.6634 & 0.6204 & 0.6516 & 0.6101 & 0.5332 & 0.4594 & 0.3712 \\ 
  2 & 0.4831 & 0.5289 & 0.5987 & 0.5828 & 0.5212 & 0.4540 & 0.3690 \\ 
  w & 0.2066 & 0.2113 & 0.2463 & 0.2918 & 0.2918 & 0.2617 & 0.2317 \\ 
  m & 0.5864 & 0.3543 & 0.6548 & 0.5325 & 0.5496 & 0.4013 & 0.3769 \\ 
  n & 1.2788 & 0.6335 & 0.2680 & 0.7180 & 0.9890 & 0.3297 & 0.3145 \\ 
   \hline
\end{tabular} \label{MSE2.5}
\end{table}

\begin{table}[tbp]
\centering
\caption{MSE of asymptotic variance estimators; $a_n = b_n^{3}$}
\begin{tabular}{rrrrrrrr}
  \hline
$*$ & $b_n = 2^4$ & $b_n = 2^5$ & $b_n = 2^6$ & $b_n = 2^7$ & $b_n = 2^8$ & $b_n = 2^9$ & $b_n = 2^{10}$ \\ 
  \hline
1 & 0.5676 & 0.4363 & 0.2485 & 0.1052 & 0.0364 & 0.0113 & 0.0038 \\ 
  2 & 0.4185 & 0.3663 & 0.2259 & 0.0994 & 0.0351 & 0.0110 & 0.0037 \\ 
  w & 0.1854 & 0.1289 & 0.0736 & 0.0344 & 0.0106 & 0.0025 & 0.0013 \\ 
  m & 0.5149 & 0.2247 & 0.2333 & 0.0798 & 0.0348 & 0.0069 & 0.0034 \\ 
  n & 1.1188 & 0.4379 & 0.0818 & 0.1119 & 0.0885 & 0.0037 & 0.0018 \\ 
   \hline
\end{tabular} \label{MSE3.0}
\end{table}

\begin{table}[tbp]
\centering
\caption{MSE of asymptotic variance estimators; $a_n = b_n^{3.5}$}
\begin{tabular}{rrrrrrrr}
  \hline
$*$ & $b_n = 2^4$ & $b_n = 2^5$ & $b_n = 2^6$ & $b_n = 2^7$ & $b_n = 2^8$ & $b_n = 2^9$ & $b_n = 2^{10}$ \\ 
  \hline
1 & 0.3284 & 0.1144 & 0.0423 & 0.0135 & 0.0058 & 0.0027 & 0.0015 \\ 
  2 & 0.2468 & 0.0948 & 0.0391 & 0.0128 & 0.0057 & 0.0027 & 0.0015 \\ 
  w & 0.1003 & 0.0288 & 0.0132 & 0.0082 & 0.0059 & 0.0043 & 0.0027 \\ 
  m & 0.2784 & 0.0471 & 0.0313 & 0.0086 & 0.0041 & 0.0017 & 0.0009 \\ 
  n & 0.6123 & 0.1014 & 0.0142 & 0.0125 & 0.0147 & 0.0025 & 0.0010 \\ 
   \hline
\end{tabular} \label{MSE3.5}
\end{table}


\begin{table}[tbp]
\centering
\caption{$b_n$ $\times$ MSE of asymptotic variance estimators; $a_n = b_n^{2}$}
\begin{tabular}{rrrrrrrr}
  \hline
$*$ & $b_n = 2^4$ & $b_n = 2^5$ & $b_n = 2^6$ & $b_n = 2^7$ & $b_n = 2^8$ & $b_n = 2^9$ & $b_n = 2^{10}$ \\ 
  \hline
1 & 10.4229 & 21.6803 & 44.3604 & 87.4912 & 176.4313 & 354.4660 & 707.3375 \\ 
  2 & 7.6130 & 18.6801 & 41.1353 & 83.6707 & 172.7787 & 350.4819 & 702.9933 \\ 
  w & 3.3640 & 7.2915 & 17.8084 & 45.4855 & 101.8822 & 215.6920 & 477.9334 \\ 
  m & 9.7003 & 12.7150 & 47.6132 & 82.8206 & 189.1759 & 322.4385 & 736.0519 \\ 
  n & 21.3576 & 23.5368 & 20.4572 & 114.4274 & 352.5406 & 290.7497 & 678.8655 \\ 
   \hline
\end{tabular} \label{bnMSE2.0}
\end{table}

\begin{table}[tbp]
\centering
\caption{$b_n$ $\times$ MSE of asymptotic variance estimators; $a_n = b_n^{2.5}$}
\begin{tabular}{rrrrrrrr}
  \hline
$*$ & $b_n = 2^4$ & $b_n = 2^5$ & $b_n = 2^6$ & $b_n = 2^7$ & $b_n = 2^8$ & $b_n = 2^9$ & $b_n = 2^{10}$ \\ 
  \hline
1 & 10.6144 & 19.8530 & 41.7007 & 78.0909 & 136.5013 & 235.2049 & 380.0606 \\ 
  2 & 7.7301 & 16.9243 & 38.3189 & 74.5957 & 133.4145 & 232.4373 & 377.8717 \\ 
  w & 3.3050 & 6.7621 & 15.7652 & 37.3484 & 74.6952 & 134.0158 & 237.2865 \\ 
  m & 9.3827 & 11.3368 & 41.9059 & 68.1650 & 140.6983 & 205.4797 & 385.9489 \\ 
  n & 20.4606 & 20.2736 & 17.1543 & 91.9000 & 253.1752 & 168.8176 & 322.0486 \\ 
   \hline
\end{tabular} \label{bnMSE2.5}
\end{table}

\begin{table}[tbp]
\centering
\caption{$b_n$ $\times$ MSE of asymptotic variance estimators; $a_n = b_n^{3}$}
\begin{tabular}{rrrrrrrr}
  \hline
$*$ & $b_n = 2^4$ & $b_n = 2^5$ & $b_n = 2^6$ & $b_n = 2^7$ & $b_n = 2^8$ & $b_n = 2^9$ & $b_n = 2^{10}$ \\ 
  \hline
1 & 9.0823 & 13.9603 & 15.9025 & 13.4609 & 9.3177 & 5.7789 & 3.8676 \\ 
  2 & 6.6952 & 11.7201 & 14.4583 & 12.7193 & 8.9962 & 5.6229 & 3.7991 \\ 
  w & 2.9661 & 4.1263 & 4.7102 & 4.4094 & 2.7073 & 1.2673 & 1.3712 \\ 
  m & 8.2390 & 7.1898 & 14.9288 & 10.2106 & 8.9112 & 3.5153 & 3.4378 \\ 
  n & 17.9001 & 14.0135 & 5.2371 & 14.3203 & 22.6665 & 1.9175 & 1.7924 \\ 
   \hline
\end{tabular} \label{bnMSE3.0}
\end{table}

\begin{table}[tbp]
\centering
\caption{$b_n$ $\times$ MSE of asymptotic variance estimators; $a_n = b_n^{3.5}$}
\begin{tabular}{rrrrrrrr}
  \hline
$*$ & $b_n = 2^4$ & $b_n = 2^5$ & $b_n = 2^6$ & $b_n = 2^7$ & $b_n = 2^8$ & $b_n = 2^9$ & $b_n = 2^{10}$ \\ 
  \hline
1 & 5.2550 & 3.6621 & 2.7055 & 1.7321 & 1.4969 & 1.3790 & 1.4981 \\ 
  2 & 3.9496 & 3.0341 & 2.4992 & 1.6330 & 1.4694 & 1.3785 & 1.4858 \\ 
  w & 1.6055 & 0.9202 & 0.8457 & 1.0472 & 1.4993 & 2.1785 & 2.7900 \\ 
  m & 4.4538 & 1.5082 & 2.0064 & 1.0981 & 1.0475 & 0.8577 & 0.8826 \\ 
  n & 9.7961 & 3.2457 & 0.9092 & 1.5997 & 3.7559 & 1.2657 & 1.0496 \\ 
   \hline
\end{tabular} \label{bnMSE3.5}
\end{table}

Table \ref{MSE2.0}, \ref{MSE2.5}, \ref{MSE3.0} and \ref{MSE3.5} 
display the MSE of the asymptotic variance estimators 
in the case that $a_n = b_n^r$ ($r = 2, 2.5, 3, 3.5$), respectively. 
Table \ref{MSE2.0} and \ref{MSE2.5} 
show that the MSE of the estimators do not converge to $0$. 
Table \ref{MSE3.0} and \ref{MSE3.5} show that they converge to $0$. 
These support Corollary \ref{avar} (a) and Corollary \ref{avar_k} (a). 

Table \ref{bnMSE2.0}, \ref{bnMSE2.5}, \ref{bnMSE3.0} and \ref{bnMSE3.5} 
display $b_n$ times the MSE of the asymptotic variance estimators 
in the case that $a_n = b_n^r$ ($r = 2, 2.5, 3, 3.5$), respectively. 
Table \ref{bnMSE2.0} and \ref{bnMSE2.5} 
show that $b_n$ times the MSE of the estimators do not converge to $0$. 
Table \ref{bnMSE3.0} and \ref{bnMSE3.5} 
show that $b_n$ times the MSE of $\hat{\Xi}_{n,*}$ ($* = 1, 2, m, n$) converge to $0$, 
and $b_n$ times the MSE of $\hat{\Xi}_{n,w}$ does not. 
These support Corollary \ref{avar} (c) and Corollary \ref{avar_k} (c).

\section*{Acknowledgment} 
The author would like to thank his supervisor Prof. Nakahiro Yoshida who provided valuable suggestions and comments. 


%
%

\end{document}